\makeatletter

\newdimen\paperwidth
\newdimen\paperheight

\def\papersize#1#2{\let\p@persize\relax\paperwidth#1\paperheight#2}

\def\Afour{\papersize{210truemm}{297truemm}}

\let\p@persize\Afour

\let\onesidestyle\@twosidefalse
\let\twosidestyle\@twosidetrue

\def\margins{\@ifnextchar[{\@margins}{\@margins[\z@]}}

\def\@margins[#1]#2#3{
  \p@persize\dimen0 #3\dimen0 .5\dimen0\normalsize%
  \oddsidemargin-1truein\advance\oddsidemargin#2%
  \evensidemargin-1truein\advance\evensidemargin#2%
  \topmargin-1truein\advance\topmargin\dimen0\headsep\dimen0\footskip\dimen0%
  \textwidth\paperwidth\advance\textwidth-#2\advance\textwidth-#2%
  \textheight\paperheight\advance\textheight-#3\advance\textheight-#3%
  \headheight\baselineskip\advance\topmargin-.5\baselineskip%
  \advance\headsep-.5\baselineskip%
  \footheight\baselineskip
  \advance\textwidth-#1\advance\oddsidemargin#1
  \if@twoside\def\@themargin%
    {\ifodd\count\z@\oddsidemargin\else\evensidemargin\fi}\fi}

\def\headlinesep#1{\advance\topmargin\headsep\advance\topmargin -#1
  \advance\topmargin.5\baselineskip\headsep #1\advance\headsep-.5\baselineskip}

\def\headline{\if@twoside\let\n@xt\h@dlin@\else\let\n@xt\h@@dlin@\fi\n@xt}
  
\def\h@dlin@#1#2{%
  \def\@oddhead{%
    {{\leftskip\z@\rightskip\z@\noindent\normalsize#1}}}
  \def\@evenhead{%
    {{\leftskip\z@\rightskip\z@\noindent\normalsize#2}}}}

\def\h@@dlin@#1{%
  \def\@oddhead{{{\leftskip\z@\rightskip\z@\noindent\normalsize#1}}}}

\def\footline{\if@twoside\let\n@xt\f@tlin@\else\let\n@xt\f@@tlin@\fi\n@xt}

\def\f@tlin@#1#2{%
  \def\@oddfoot{%
    {{\leftskip\z@\rightskip\z@\noindent\normalsize#1}}}
  \def\@evenfoot{%
    {{\leftskip\z@\rightskip\z@\noindent\normalsize#2}}}}

\def\f@@tlin@#1{%
  \def\@oddfoot{{{\leftskip\z@\rightskip\z@\noindent\normalsize#1}}}}

\def\normalpage{\global\@specialpagefalse}

\makeatother

\makeatletter

\def\ft{\@ifnextchar[{\ft@s}{\ft@}}
\def\ft@{\ft@@@s[\f@size]}
\def\ft@s[{\@ifnextchar{a}{\ft@sz[}{\ft@@s[}}
\def\ft@@s[{\@ifnextchar{s}{\ft@sz[}{\ft@@@s[}}
\def\ft@@@s[#1]{\ft@sz[at #1pt]}
\def\ft@sz[#1]#2{\font\fonttemp=#2 #1\fonttemp\ignorespaces}

\makeatother

\makeatletter

\input{epsf.sty}

\makeatother

\def\smallcircc{\mathop{\mkern3.5mu\hbox{\raise.58ex\hbox{\ft{lcircle10}a}}}}
\def\varemptyset{{\hbox{\raise.21ex\hbox{$\not$}}\mkern.15mu\mathrm{O}\mkern.15mu}}

  \let\epsilon\varepsilon
      \let\theta\vartheta
          \let\phi\varphi
   \let\emptyset\varemptyset

\documentstyle[12pt]{article}

\makeatletter

\let\Larg@\Large
\let\hug@\huge

\def\usepackage#1{\input{#1.sty}}

\input{geom.sty}

\def\r@adlabel#1#2{\global\@namedef{#1@\the\@key}{#2}}

\let\Large\Larg@
\let\huge\hug@

\def\smallskip{\vskip\smallskipamount}
\def\medskip{\vskip\medskipamount}
\def\bigskip{\vskip\bigskipamount}

\def\mytrivlist{\parsep\parskip\@nmbrlistfalse
  \my@trivlist \labelwidth\z@ \leftmargin\z@
  \itemindent\z@ \def\makelabel##1{##1}}

\def\my@trivlist{\global\@newlisttrue \@outerparskip\parskip}

\def\end#1{\csname end#1\endcsname\@checkend{#1}%
  \expandafter\endgroup\if@endpe\@doendpe\fi
  \if@ignore \global\@ignorefalse \ignorespaces\fi}

\def\put{\@ifnextchar[{\@put}{\@@rput[\z@,\z@][r]}}
\def\@put[#1]{\@ifnextchar[{\@@put[#1]}{\@@@@@put[#1]}}
\def\@@put[#1][{\@ifnextchar{l}{\@@lput[#1][}{\@@@put[#1][}}
\def\@@@put[#1][{\@ifnextchar{c}{\@@cput[#1][}{\@@@@put[#1][}}
\def\@@@@put[#1][{\@ifnextchar{r}{\@@rput[#1][}{\relax}}
\def\@@@@@put[{\@ifnextchar{l}{\@@lput[\z@,\z@][}{\@@@@@@put[}}
\def\@@@@@@put[{\@ifnextchar{c}{\@@cput[\z@,\z@][}{\@@@@@@@put[}}
\def\@@@@@@@put[{\@ifnextchar{r}{\@@rput[\z@,\z@][}{\@@@@@@@@put[}}
\def\@@@@@@@@put[#1]{\@@rput[#1][r]}

\let\hm@d@\leavevmode

\long\def\@@lput[#1,#2][l]#3{\setbox0\hbox{#3}\hm@d@\raise#2\hbox to\z@{\dimen0 #1%
  \advance\dimen0-\wd0\kern\dimen0\dp0\z@\ht0\z@\wd0\z@\box0\hss}\ignorespaces}
\long\def\@@cput[#1,#2][c]#3{\setbox0\hbox{#3}\hm@d@\raise#2\hbox to\z@{\dimen0 #1%
  \advance\dimen0-.5\wd0\kern\dimen0\dp0\z@\ht0\z@\wd0\z@\box0\hss}\ignorespaces}
\long\def\@@rput[#1,#2][r]#3{\setbox0\hbox{\kern#1\raise#2\hbox{#3}}%
  \dp0\z@\ht0\z@\wd0\z@\hm@d@\box0\ignorespaces}

\def\flbox{\@ifnextchar[{\@flbox}{\@@rflbox[\z@,\z@][r]}}
\def\@flbox[#1]{\@ifnextchar[{\@@flbox[#1]}{\@@@@@flbox[#1]}}
\def\@@flbox[#1][{\@ifnextchar{l}{\@@lflbox[#1][}{\@@@flbox[#1][}}
\def\@@@flbox[#1][{\@ifnextchar{c}{\@@cflbox[#1][}{\@@@@flbox[#1][}}
\def\@@@@flbox[#1][{\@ifnextchar{r}{\@@rflbox[#1][}{\relax}}
\def\@@@@@flbox[{\@ifnextchar{l}{\@@lflbox[\z@,\z@][}{\@@@@@@flbox[}}
\def\@@@@@@flbox[{\@ifnextchar{c}{\@@cflbox[\z@,\z@][}{\@@@@@@@flbox[}}
\def\@@@@@@@flbox[{\@ifnextchar{r}{\@@rflbox[\z@,\z@][}{\@@@@@@@@flbox[}}
\def\@@@@@@@@flbox[#1]{\@@rflbox[#1][r]}
\long\def\@@lflbox[#1,#2][l]#3{\@@lput[#1,#2][l]{%
  \vtop{\leftskip\z@\parindent\z@\raggedleft\hm@d@#3}}}
\long\def\@@cflbox[#1,#2][c]#3{\@@cput[#1,#2][c]{%
  \vtop{\leftskip\z@\parindent\z@\raggedcenter\hm@d@#3}}}
\long\def\@@rflbox[#1,#2][r]#3{\@@rput[#1,#2][r]{%
  \vtop{\leftskip\z@\parindent\z@\raggedright\hm@d@#3}}}

\def\maketitle{\par
 \begingroup
 \def\thefootnote{\fnsymbol{footnote}}
 \def\@makefnmark{\hbox to 0pt{$^{\@thefnmark}$\hss}} 
 \if@twocolumn 
 \twocolumn[\@maketitle] 
 \else 
 \global\@topnum\z@ \@maketitle \fi\thispagestyle{plain}\@thanks
 \endgroup
 \setcounter{footnote}{0}
 \let\maketitle\relax
 \let\@maketitle\relax
 \gdef\@thanks{}\gdef\@author{}\gdef\@title{}\let\thanks\relax}

\def\@maketitle{ 
 \null
 \vskip 2em \begin{center}
 {\LARGE \@title \par} \vskip 1.5em {\large \lineskip .5em
\begin{tabular}[t]{c}\@author 
 \end{tabular}\par} 
 \vskip 1em {\large \@date} \end{center}
 \par
 \vskip 1.5em}
 
\def\partbeforeskip#1{\def\p@rtbeforeskip{#1}}
\def\partstyle#1{\def\p@rtstyl@{#1}}
\def\partdot#1{\def\partd@t{#1}}
\def\partafterskip#1{\def\p@rtafterskip{#1}}
\def\partintrostyle#1{\def\partintr@styl@{#1}}
\def\partintrodot#1{\def\partintr@dot{#1}}
\long\def\partintrosep#1{\long\def\partintr@sep{#1}}
\def\partnewpagetrue{\def\p@rtnewp@ge{\newpage}}
\def\partnewpagefalse{\long\def\p@rtnewp@ge{\par}}

\partbeforeskip{4ex}
\partstyle{\centering\Large\bf}
\partdot{}
\partafterskip{3ex}
\partintrostyle{\large}
\partintrodot{}
\partintrosep{\par}
\partnewpagefalse

\def\partname{Part}
\def\part{\p@rtnewp@ge\addvspace\p@rtbeforeskip\@afterindentfalse\secdef\@part\@spart}

\def\@part[#1]#2{\ifnum \c@secnumdepth >-1\relax  
        \refstepcounter{part}                     
        \def\@tempa{\addcontentsline{toc}{part}}  %
        \expandafter\@tempa\expandafter{\thepart  
          \hspace{1em}#1}\else                    
        \addcontentsline{toc}{part}{#1}\fi        
   {\p@rtstyl@                       
    \ifnum \c@secnumdepth >-1\relax        
      {\partintr@styl@\partname\ \thepart  
       \partintr@dot}\partintr@sep\nobreak 
    \fi                                    
    #2\partd@t\markboth{}{}\par}
    \nobreak                       
    \vskip\p@rtafterskip           
   \@afterheading                  
    }                              

\def\@spart#1{{\p@rtcentering\p@rtstyl@                      
    #1\partd@t\par}                 
    \nobreak                        
    \vskip\p@rtafterskip            
    \@afterheading                  
  }                                 

\newif\ifsection@ftind
\newif\ifsection@ftpar

\def\sectionbeforeskip#1{\def\s@ctbeforeskip{#1}}
\def\sectionstyle#1{\def\s@ctstyl@{#1}}
\def\sectiondot#1{\def\sectiond@t{#1}}
\def\sectionafterskip#1{\def\s@ctafterskip{#1}}
\def\sectionintrostyle#1{\def\sectionintr@styl@{#1}}
\def\sectionintro#1{\def\sectionintr@{#1}}
\def\sectionintrodot#1{\def\sectionintr@dot{#1}}
\def\sectionintrosep#1{\def\sectionintr@sep{#1}}
\def\sectionindenttrue{\def\s@ctind{\parindent}}
\def\sectionindentfalse{\def\s@ctind{\z@}}
\def\sectionafterindenttrue{\section@ftindtrue}
\def\sectionafterindentfalse{\section@ftindfalse}
\def\sectionafternewlinetrue{\section@ftpartrue}
\def\sectionafternewlinefalse{\section@ftparfalse}

\newif\ifsubsection@ftind
\newif\ifsubsection@ftpar

\def\subsectionbeforeskip#1{\def\ss@ctbeforeskip{#1}}
\def\subsectionstyle#1{\def\ss@ctstyl@{#1}}
\def\subsectiondot#1{\def\subsectiond@t{#1}}
\def\subsectionafterskip#1{\def\ss@ctafterskip{#1}}
\def\subsectionintrostyle#1{\def\subsectionintr@styl@{#1}}
\def\subsectionintro#1{\def\subsectionintr@{#1}}
\def\subsectionintrodot#1{\def\subsectionintr@dot{#1}}
\def\subsectionintrosep#1{\def\subsectionintr@sep{#1}}
\def\subsectionindenttrue{\def\ss@ctind{\parindent}}
\def\subsectionindentfalse{\def\ss@ctind{\z@}}
\def\subsectionafterindenttrue{\subsection@ftindtrue}
\def\subsectionafterindentfalse{\subsection@ftindfalse}
\def\subsectionafternewlinetrue{\subsection@ftpartrue}
\def\subsectionafternewlinefalse{\subsection@ftparfalse}

\newif\ifsubsubsection@ftind
\newif\ifsubsubsection@ftpar

\def\subsubsectionbeforeskip#1{\def\sss@ctbeforeskip{#1}}
\def\subsubsectionstyle#1{\def\sss@ctstyl@{#1}}
\def\subsubsectiondot#1{\def\subsubsectiond@t{#1}}
\def\subsubsectionafterskip#1{\def\sss@ctafterskip{#1}}
\def\subsubsectionintrostyle#1{\def\subsubsectionintr@styl@{#1}}
\def\subsubsectionintro#1{\def\subsubsectionintr@{#1}}
\def\subsubsectionintrodot#1{\def\subsubsectionintr@dot{#1}}
\def\subsubsectionintrosep#1{\def\subsubsectionintr@sep{#1}}
\def\subsubsectionindenttrue{\def\sss@ctind{\parindent}}
\def\subsubsectionindentfalse{\def\sss@ctind{\z@}}
\def\subsubsectionafterindenttrue{\subsubsection@ftindtrue}
\def\subsubsectionafterindentfalse{\subsubsection@ftindfalse}
\def\subsubsectionafternewlinetrue{\subsubsection@ftpartrue}
\def\subsubsectionafternewlinefalse{\subsubsection@ftparfalse}

\newif\ifparagraph@ftind
\newif\ifparagraph@ftpar

\def\paragraphbeforeskip#1{\def\p@rbeforeskip{#1}}
\def\paragraphstyle#1{\def\p@rstyl@{#1}}
\def\paragraphdot#1{\def\paragraphd@t{#1}}
\def\paragraphafterskip#1{\def\p@rafterskip{#1}}
\def\paragraphintrostyle#1{\def\paragraphintr@styl@{#1}}
\def\paragraphintro#1{\def\paragraphintr@{#1}}
\def\paragraphintrodot#1{\def\paragraphintr@dot{#1}}
\def\paragraphintrosep#1{\def\paragraphintr@sep{#1}}
\def\paragraphindenttrue{\def\p@rind{\parindent}}
\def\paragraphindentfalse{\def\p@rind{\z@}}
\def\paragraphafterindenttrue{\paragraph@ftindtrue}
\def\paragraphafterindentfalse{\paragraph@ftindfalse}
\def\paragraphafternewlinetrue{\paragraph@ftpartrue}
\def\paragraphafternewlinefalse{\paragraph@ftparfalse}

\newif\ifsubparagraph@ftind
\newif\ifsubparagraph@ftpar

\def\subparagraphbeforeskip#1{\def\sp@rbeforeskip{#1}}
\def\subparagraphstyle#1{\def\sp@rstyl@{#1}}
\def\subparagraphdot#1{\def\subparagraphd@t{#1}}
\def\subparagraphafterskip#1{\def\sp@rafterskip{#1}}
\def\subparagraphintrostyle#1{\def\subparagraphintr@styl@{#1}}
\def\subparagraphintro#1{\def\subparagraphintr@{#1}}
\def\subparagraphintrodot#1{\def\subparagraphintr@dot{#1}}
\def\subparagraphintrosep#1{\def\subparagraphintr@sep{#1}}
\def\subparagraphindenttrue{\def\sp@rind{\parindent}}
\def\subparagraphindentfalse{\def\sp@rind{\z@}}
\def\subparagraphafterindenttrue{\subparagraph@ftindtrue}
\def\subparagraphafterindentfalse{\subparagraph@ftindfalse}
\def\subparagraphafternewlinetrue{\subparagraph@ftpartrue}
\def\subparagraphafternewlinefalse{\subparagraph@ftparfalse}

\sectionbeforeskip{\bigskipamount}
\sectionstyle{\large\bf}
\sectiondot{}
\sectionafterskip{.5\bigskipamount}
\sectionintrostyle{}
\sectionintro{}
\sectionintrodot{.}
\sectionintrosep{1.25ex}
\sectionindentfalse
\sectionafterindenttrue
\sectionafternewlinetrue

\subsectionbeforeskip{.8\bigskipamount}
\subsectionstyle{\normalsize\bf}
\subsectiondot{}
\subsectionafterskip{.4\bigskipamount}
\subsectionintrostyle{}
\subsectionintro{}
\subsectionintrodot{.}
\subsectionintrosep{1.25ex}
\subsectionindentfalse
\subsectionafterindenttrue
\subsectionafternewlinetrue

\subsubsectionbeforeskip{.6\bigskipamount}
\subsubsectionstyle{\normalsize\bf}
\subsubsectiondot{}
\subsubsectionafterskip{.3\bigskipamount}
\subsubsectionintrostyle{}
\subsubsectionintro{}
\subsubsectionintrodot{.}
\subsubsectionintrosep{1.25ex}
\subsubsectionindentfalse
\subsubsectionafterindenttrue
\subsubsectionafternewlinetrue

\paragraphbeforeskip{.5\bigskipamount}
\paragraphstyle{\normalsize\bf}
\paragraphdot{.}
\paragraphafterskip{1.25ex}
\paragraphintrostyle{}
\paragraphintro{}
\paragraphintrodot{.}
\paragraphintrosep{1.25ex}
\paragraphindentfalse
\paragraphafterindenttrue
\paragraphafternewlinefalse

\subparagraphbeforeskip{.5\bigskipamount}
\subparagraphstyle{\normalsize\bf}
\subparagraphdot{.}
\subparagraphafterskip{1.25ex}
\subparagraphintrostyle{}
\subparagraphintro{}
\subparagraphintrodot{.}
\subparagraphintrosep{1.25ex}
\subparagraphindenttrue
\subparagraphafterindenttrue
\subparagraphafternewlinefalse

\let\@partoken\par
\long\def\@@gobble#1{}
\def\ignorepar{\@ifnextchar\@partoken{\expandafter\ignorepar\@@gobble}{\ignorespaces}}

\def\@startsection#1#2#3#4#5#6{
   \@tempskipa #4\relax
   \csname if#1@ftind\endcsname\@afterindenttrue\else\@afterindentfalse\fi
   \advance\@tempskipa by\presection
   \if@nobreak \everypar{}\else
     \addpenalty{\@secpenalty}\addvspace{\@tempskipa}%
     \allowbreak\vskip -\presection \fi \@ifstar
     {\@ssect{#1}{#2}{#3}{#4}{#5}{#6}}{\@dblarg{\@sect{#1}{#2}{#3}{#4}{#5}{#6}}}}

\def\@sect#1#2#3#4#5#6[#7]#8{\def\object@type{#1}%
   \ifnum #2>\c@secnumdepth\def\@svsec{}\def\@tempb{}%
      \else\refstepcounter{#1}\def\@svsec{{\csname #1intr@styl@\endcsname%
        {\csname #1intr@\endcsname}\csname the#1\endcsname%
        \csname #1intr@dot\endcsname\kern\csname #1intr@sep\endcsname}}%
        \edef\@tempb{\noexpand\numberline{\csname the#1\endcsname}}\fi%
   \def\@tempa{\addcontentsline{toc}{#1}}%
   \csname if#1@ftpar\endcsname%
      \begingroup #6\relax%
        \@hangfrom{\hskip #3\relax\@svsec}{\interlinepenalty \@M{#8}%
        \csname #1d@t\endcsname\par}%
      \endgroup%
      \csname #1mark\endcsname{#7}%
      \expandafter\@tempa\expandafter{\@tempb #7}%
      \ifautolabel\label*{#8}\fi%
   \else%
      \def\@svsechd{#6\hskip #3\relax%
         \@svsec{#8}\csname #1mark\endcsname{#7}%
         \expandafter\@tempa\expandafter{\@tempb #7}%
         \ifautolabel\label*{#8}\fi}\fi%
   \@xsect{#1}{#5}\ignorepar}

\def\@ssect#1#2#3#4#5#6#7{%
   \ifnum #2>\c@secnumdepth\def\@tempb{}\else \def\@tempb{\numberline{}}\fi%
     \def\@tempa{\addcontentsline{toc}{s#1}}%
     \csname if#1@ftpar\endcsname
        \begingroup #6\relax
           \@hangfrom{\hskip #3}{\interlinepenalty \@M{#7}%
           \csname #1d@t\endcsname\par}%
        \endgroup
        \csname s#1mark\endcsname{#7}%
        \ifstarredcontents\expandafter\@tempa\expandafter{\@tempb #7}\fi%
        \ifautolabel\label*{#7}\fi%
     \else%
        \def\@svsechd{#6\hskip #3\relax{#7}\csname s#1mark\endcsname{#7}%
        \ifautolabel\label*{#7}\fi}\fi
   \@xsect{#1}{#5}\ignorepar}

\def\@xsect#1#2{
   \csname if#1@ftpar\endcsname 
       \par \nobreak \vskip #2\relax \@afterheading
    \else \global\@nobreakfalse \global\@noskipsectrue
       \everypar{\if@noskipsec \global\@noskipsecfalse
                   \clubpenalty\@M \hskip -\parindent
                   \begingroup \@svsechd \endgroup \unskip
                   \hskip #2\relax  
                  \else \clubpenalty \@clubpenalty
                    \everypar{}\fi}\fi\ignorespaces}

\def\section{\@startsection{section}{1}{\s@ctind}
  {\s@ctbeforeskip}{\s@ctafterskip}{\s@ctstyl@}}
\def\subsection{\@startsection{subsection}{2}{\ss@ctind}
  {\ss@ctbeforeskip}{\ss@ctafterskip}{\ss@ctstyl@}}
\def\subsubsection{\@startsection{subsubsection}{3}{\sss@ctind}
  {\sss@ctbeforeskip}{\sss@ctafterskip}{\sss@ctstyl@}}
\def\paragraph{\@startsection{paragraph}{4}{\p@rind}
  {\p@rbeforeskip}{\p@rafterskip}{\p@rstyl@}}
\def\subparagraph{\@startsection{subparagraph}{4}{\sp@rind}
  {\sp@rbeforeskip}{\sp@rafterskip}{\sp@rstyl@}}

\def\statementabove#1{\def\th@bove{#1}}
\def\statementstyle#1{\def\thstyl@{#1}}
\def\statementbelow#1{\def\thb@low{#1}}
\def\statementindentfalse{\let\thind@nt\relax}
\def\statementindenttrue{\let\thind@nt\indent}

\def\statementintrostyle#1{\def\thintr@style{#1}}
\def\statementintrodot#1{\def\thintr@dot{#1}}
\def\statementintrosep#1{\def\thintr@sep{#1}}
\def\statementintrobrackets#1#2{\def\thintr@left{#1}\def\thintr@right{#2}}

\statementabove{\medskip}
\statementstyle{\sl}
\statementbelow{\medskip}
\statementindenttrue

\statementintrostyle{\normalshape\bf}
\statementintrodot{.}
\statementintrosep{\kern1.25ex}
\statementintrobrackets{(}{)}

\def\@thskip{\dimen0\lastskip\vskip-\dimen0%
  \th@bove\dimen1\lastskip\vskip-\dimen1%
  \ifdim\dimen0>\dimen1\else\dimen0\dimen1\fi\vskip\dimen0}

\long\def\@@newtheorem#1#2#3{%
  \newenvironment{#3}%
    {\def\object@type{#3}\par\@thskip%
     \@ifnextchar[{\@enva{#3}{\thstyl@#1{#2}}}{\@envb{#3}{\thstyl@#1{#2}}}}%
    {\end{#3@}}%
  \@ifnextchar[{\@nothm{#3}}{\@nnthm{#3}}}

\def\@nothm#1[#2]#3{%
  \@ifundefined{c@#2}{\@latexerr{No theorem environment `#2' defined}\@eha}%
  {\expandafter\@ifdefinable\csname #1@\endcsname
  {\global\@namedef{the#1}{\@nameuse{the#2}}%
   \global\@namedef{c@#1}{\@nameuse{c@#2}}
   \global\@namedef{p@#1}{\@nameuse{p@#2}}
   \global\@namedef{#1@}{\@nnnthm{#2}{#3}}%
   \global\@namedef{end#1@}{\@endtheorem}}}}

\def\@nnnthm#1#2{\refstepcounter
    {#1}\@ifnextchar[{\@ynnnthm{#1}{#2}}{\@xnnnthm{#1}{#2}}}
 
\def\@xnnnthm#1#2{\@begintheorem{#2}{\csname the#1\endcsname}\ignorespaces}
\def\@ynnnthm#1#2[#3]{\@opargbegintheorem{#2}{\csname
       the#1@\endcsname}{#3}\ignorespaces}

\def\renewtheorem{\@ifnextchar[{\@renewtheorem}{\@renewtheorem[{}{}]}}

\long\def\@renewtheorem[#1]{\@@renewtheorem#1}

\long\def\@@renewtheorem#1#2#3{%
  \expandafter\let\csname#3@\endcsname\undefined
  \renewenvironment{#3}%
    {\def\object@type{#3}\par\@thskip%
     \@ifnextchar[{\@enva{#3}{\thstyl@#1{#2}}}{\@envb{#3}{\thstyl@#1{#2}}}}%
    {\end{#3@}}%
  \@ifnextchar[{\@nothm{#3}}{\@nnthm{#3}}}

\def\@begintheorem#1#2{\@opargbegintheorem{#1}{#2}{}}

\def\@opargbegintheorem#1#2#3{%
        \def\@tempx{#1}%
        \expandafter\let\expandafter\@tempy#2
        \def\@tempz{#3}%
        \mytrivlist\item[\thind@nt\hskip\labelsep%
        {\thintr@style%
          #1\if\@tempx\@empty\else\if\@tempy\relax\else\kern1ex\fi\fi#2%
          \ifx\@tempz\@empty%
            \if\@tempx\@empty\if\@tempy\relax%
            \else\thintr@dot\thintr@sep\fi\else\thintr@dot\thintr@sep\fi%
            \else%
            \if\@tempx\@empty\if\@tempy\relax\else\kern1ex\fi\else\kern1ex\fi%
           \thintr@left{#3}\thintr@right\thintr@dot\thintr@sep\fi}%
            \hskip-\labelsep]%
        \ifautolabel\label*{#3}\fi}

\def\@endtheorem{\strut\endtrivlist\thb@low}

\def\proofabove#1{\def\pf@bove{#1}}
\def\proofstyle#1{\def\pfstyl@{#1}}
\def\proofbelow#1{\def\pfb@low{#1}}
\def\proofindentfalse{\let\pfind@nt\relax}
\def\proofindenttrue{\let\pfind@nt\indent}

\def\proofintrostyle#1{\def\pfintr@style{#1}}
\def\proofintrodot#1{\def\pfintr@dot{#1}}
\def\proofintrosep#1{\def\pfintr@sep{#1}}
\def\proofintrobrackets#1#2{\def\pfintr@left{#1}\def\pfintr@right{#2}}

\proofabove{\medskip}
\proofstyle{}
\proofbelow{\medskip}
\proofindenttrue

\proofintrostyle{\sl}
\proofintrodot{.}
\proofintrosep{\kern1.25ex}
\proofintrobrackets{of\kern1ex}{}

\def\@pfskip{\dimen0\lastskip\vskip-\dimen0%
  \pf@bove\dimen1\lastskip\vskip-\dimen1%
  \ifdim\dimen0>\dimen1\else\dimen0\dimen1\fi\vskip\dimen0}

\renewenvironment{proof}%
  {\@pfskip\mytrivlist\item[\pfind@nt]\@ifnextchar[{\pro@f}{\pro@f[\prooftag]}}
  {\ifvoid\provedbox\else\hproved\fi\endtrivlist\pfb@low}

\def\pro@f[#1]{\setbox\provedbox\hbox{\provedboxcontents{#1}}\proofintro{#1}}

\def\proofintro#1{\expandafter\def\expandafter\@tempa\expandafter{#1}%
  {\pfintr@style{Proof\ifx\@tempa\empty\else\kern1ex\pfintr@left{#1}%
  \pfintr@right\fi}\pfintr@dot\pfintr@sep}\pfstyl@\ignorespaces}

\def\provedmark#1{\def\prm@rk{#1}}
\def\provedsep#1{\def\prs@p{#1}}

\provedmark{$\square$}
\provedsep{\kern1.25ex}

\def\provedtexttrue{\def\prb@x##1{\fbox{\small##1}}}
\def\provedtextfalse{\def\prb@x##1{\prm@rk}}
\def\provedmarkrighttrue{\let\prhf@l\hfill}
\def\provedmarkrightfalse{\let\prhf@l\relax}

\provedtextfalse
\provedmarkrightfalse

\def\provedboxcontents#1{\expandafter\def\expandafter\@tempa\expandafter{#1}%
  \ifx\@tempa\empty\prm@rk\else\prb@x{#1}\fi}

\def\proved{\ifmmode\eqno{\box\provedbox}\else\hproved\fi}

\def\hproved{\unskip\nobreak\prhf@l\penalty50\prs@p\hbox{}\nobreak\prhf@l
  \box\provedbox{\finalhyphendemerits=0\par}}

\def\captionstyle#1{\def\c@ptstyl@{#1}}
\def\captionintrostyle#1{\def\c@pintr@style{#1}}
\def\captionintrodot#1{\def\c@pintr@dot{#1}}
\def\captionintrosep#1{\def\c@pintr@sep{#1}}

\captionstyle{\small\sf}
\captionintrostyle{\bf}
\captionintrodot{.}
\captionintrosep{\hskip1.25ex}

\long\def\@makecaption#1#2{%
    \vskip\captionskip
    \setbox\@tempboxa\hbox{%
      \ifproofing\@ifundefined{the@label}{}
        {\hbox to 0pt{\vbox to 0pt{\vss\hbox{\tiny\the@label}\bigskip}\hss}}\fi
      \c@ptstyl@{\c@pintr@style #1\c@pintr@dot}\ignorespaces #2}%
    \@captionwidth=\hsize \advance\@captionwidth-2\@captionmargin
    \ifdim \wd\@tempboxa >\@captionwidth {%
        \rightskip=\@captionmargin\leftskip=\@captionmargin
        \unhbox\@tempboxa\par}%
      \else
        \hbox to\hsize{\hfil\box\@tempboxa\hfil}%
    \fi}

\def\end@Float#1{%
  \expandafter\caption\expandafter[\the@title]{%
   {\c@pintr@style%
   \ifx\the@caption\empty\ifx\the@title\empty
   \else\c@pintr@sep\fi\else\c@pintr@sep\fi
    \the@title\ifx\the@caption\empty%
     \expandafter\label\expandafter*\expandafter{\the@label}%
    \else\ifx\the@title\empty%
     \expandafter\label\expandafter*\expandafter{\the@label}%
    \else\c@pintr@dot\c@pintr@sep%
     \expandafter\label\expandafter*\expandafter{\the@label}\fi\fi}%
   \ignorespaces\the@caption}%
  \end{#1}}

\renewenvironment{Figure}{\@ifnextchar[%
  {\@myFloat{figure}}{\@myFloat{figure}[htbp]}}{\end@Float{figure}}

\def\@myFloat#1[#2]#3{%
  \begin{#1}[#2]\def\the@label{#3}}

\def\showfig{\showfigurestrue\fig}
\def\fig#1#2{\@ifnextchar[{\@fig{#1}{#2}}{\@fig{#1}{#2}[0pt]}}

\def\@fig#1#2[#3]#4#5{%
  \def\the@title{#4}\def\the@caption{#5}\centerline{\fig@{#1}{#2}}\vskip#3}

\def\fig@@#1#2{\leavevmode{\figstyl@\vrule width 0pt height 1.8ex%
 \smash{\framebox{\strut\def\@temp{#1}\ifx\@temp\@empty{ #2 }\else{ #1 }\fi}}}}
\def\fig@@@#1#2{\epsfbox{#2}}

\def\figstyle#1{\def\figstyl@{#1}}

\figstyle{\normalshape\mediumseries\normalsize}

\def\@ldshowfig#1#2{\epsfbox{#2}}
\def\@ldfig@#1#2{\leavevmode{\framebox{\figstyl@\strut{ #1 }}}}
\def\@ldshowfigurestrue{\let\fig\showfig}
\def\@ldshowfiguresfalse{\let\fig\fig@}

\def\showfiguresfalse{\let\fig@\fig@@}
\def\showfigurestrue{\let\fig@\fig@@@}

\@definecounter{bibenumi}

\def\thebibliography#1{%
 \section*{\refname}\vskip-\lastskip%
 \list{[\arabic{bibenumi}]}{\topsep0pt\settowidth\labelwidth{[#1]}%
 \leftmargin\labelwidth\advance\leftmargin\labelsep\usecounter{bibenumi}}%
 \def\newblock{\hskip .11em plus .33em minus .07em}%
 \sloppy\clubpenalty4000\widowpenalty4000\sfcode`\.=1000\relax}

\makeatother

\proofingfalse
\autolabelfalse

\showfigurestrue

\newtheorem{stat}{\statname}  \unnumbered{stat}
\newenvironment{statement}[1]{\def\statname{#1}\begin{stat}}{\end{stat}}

\newtheorem{nstat}{\nstatname}[section]

\newtheorem{lemma}[nstat]{Lemma}

\newtheorem[{\ns}{}]{remark}[nstat]{Remark}

\papersize{215.9truemm}{279.4truemm}
\margins{3.295cm}{2.12cm}

\headline{\hfill}
\footline{\small\hfill--\kern1ex\thepage\kern1ex--\hfill}
\makeatletter
  \textfloatsep\floatsep
  \intextsep\floatsep
\makeatother

\sectionbeforeskip{1.5\bigskipamount}
\sectionstyle{\centering\normalsize\bf}
\sectionafterskip{\bigskipamount}
\sectionafterindenttrue

\subsectionbeforeskip{\bigskipamount}
\subsectionstyle{\centering\normalsize\sl}
\subsectionafterskip{.5\bigskipamount}
\subsectionafterindenttrue

\paragraphbeforeskip{\bigskipamount}
\paragraphstyle{\centering\normalsize\sl}
\paragraphafterskip{.5\bigskipamount}
\paragraphafternewlinetrue
\paragraphafterindenttrue
\paragraphdot{}

\statementintrostyle{\sc}
\renewtheorem[{\ns}{}]{remark}[nstat]{Remark}
\renewtheorem[{\ns}{}]{definition}[nstat]{Definition}

\captionstyle{\small}
\captionintrostyle{\sc}

\newcommand{\id}{\mathop{\mathrm{id}}\nolimits} 
\newcommand{\Cl}{\mathop{\mathrm{Cl}}\nolimits} 
 
\newcommand{\Int}{\mathop{\mathrm{Int}}\nolimits} 
\newcommand{\Bd}{\mathop{\mathrm{Bd}}\nolimits}

\newcommand{\Fix}{\mathop{\mathrm{Fix}}\nolimits}

\let\emptyset\varemptyset

\def\(#1){$(${\sl #1}\/$)$}

\def\varemptyset{%
 {\text{\raise.21ex\hbox{$\not$}}\mkern.15mu\mathrm{O}\mkern.15mu}}
\def\widebar#1{%
 \text{\setbox0\hbox{$#1$}\dimen0\ht0\dimen0.25\dimen0
 \hbox to 0pt{$\kern\dimen0\overline{\kern-\dimen0\phantom{#1}}$\hss}}#1}

\font\ftt cmtt10 at 11pt
\font\fsc cmcsc10 at 12pt
\font\fsl cmsl12 at 12pt

\begin{document}\frenchspacing

\title{\large\bf INVOLUTIONS OF 3-DIMENSIONAL HANDLEBODIES%
       \label{Version 0.5 / \today}}
\author{
\fsc Andrea Pantaleoni\\
\fsl Dipartimento di Matematica\\[-3pt]
\fsl Universit\`a di Camerino -- Italy\\
\ftt andrea.pantaleoni@unicam.it
\and
\fsc Riccardo Piergallini\\
\fsl Dipartimento di Matematica\\[-3pt]
\fsl Universit\`a di Camerino -- Italy\\
\ftt riccardo.piergallini@unicam.it
}
\date{}

\maketitle

\begin{abstract}
\baselineskip13.5pt
\smallskip

\noindent
We study the orientation preserving involutions of the orientable 3-dimensional handlebody
$H_g$, for any genus $g$. A complete classification of such involutions is given in terms
of their fixed points.



\end{abstract}

\section*{Introduction}

Involutions of the 3-dimensional orientable handlebody $H_g$ of genus $g$ has been already
classified in \cite{Bar}, \cite{Liv}, \cite{Wal} and \cite{Prz} for $g \leq 2$. Moreover,
a classification of the orientation reversing involutions of $H_g$ was given in \cite{KMG}
(theorem 3.6).

In this paper, we complete the study of the subject, by providing a classification of the
orientation preserving involutions of $H_g$ for any $g \geq 0$. Our argument is direct and
elementary. The same result can also be derived from the general theory of actions on
handlebodies developed in \cite{CMZ}.

\smallskip

Namely, we prove the following theorem.

\begin{statement}{Theorem} 
Let $h: H_g \to H_g$ be an orientation preserving involution. If $h$ is free, then $g =
2n+1$ for some $n \geq 0$ and $h$ is equivalent to the involution $I_g$ depicted in Figure
\ref{involution1/fig}. If $h$ is not free, then there exist $n,m,l \geq 0$ with $1 \leq n
+ 2m \leq n + 2m + 2l = g + 1$, such that $h$ is equivalent to the involution $L_g^{n,m}$
depicted in Figure \ref{involution2/fig}.
\end{statement}

\begin{Figure}[b]{involution1/fig}
\fig{}{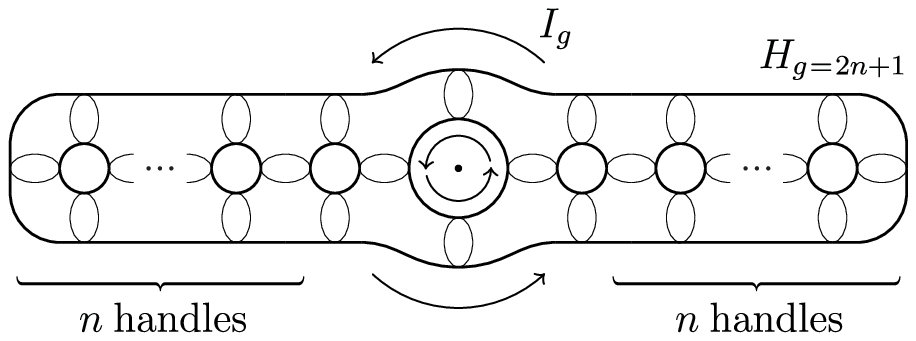}{}{The free involution $I_g$ for $g = 2n + 1$}\vskip-6pt
\end{Figure}

The free involution $I_g: H_g \to H_g$ with $g = 2n +1$ can be realized by embedding $H_g$
in $R^3$ as in Figure \ref{involution1/fig} and rotating of $\pi$ radians around the axis
orthogonal to the plane of the picture at the dot.

\begin{Figure}[htb]{involution2/fig}
\fig{}{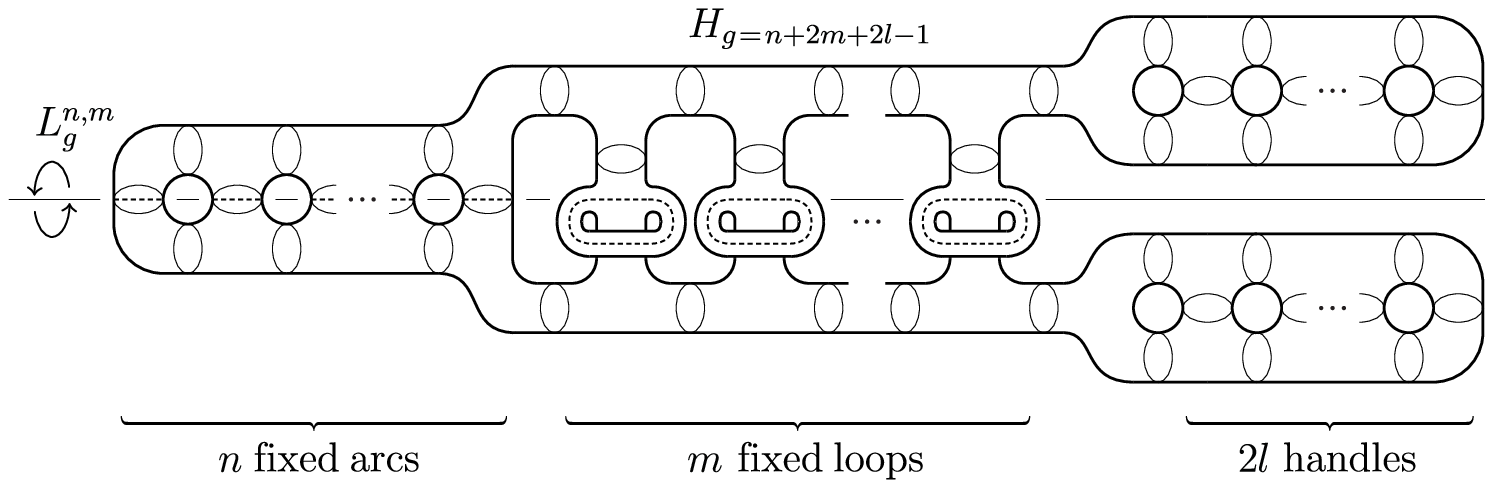}{}{The non-free involution $L^{n,m}_g$}
\end{Figure}

The description of the involution $L_g^{n,m}: H_g \to H_g$, with $g = n + 2m + 2l -1$, is
a little bit more involved. The fixed point set $\Fix L_g^{n,m}$ consists of $n$ arcs and
$m$ loops, all dashed in Figure \ref{involution2/fig}. We think of the $H_g$ as $H_{n + m
+ 2l - 1}$ with $m$ extra handles attached to it. The handlebody $H_{n + m + 2l - 1}$ is
imbedded in $R^3$ in such a way that it is symmetric with respect to the median horizontal
line and meets it in $n + m$ arcs, while the $m$ extra handles are the non-symmetric ones.
Then, the restriction of $L_g^{n,m}$ to $H_{n + m + 2l - 1}$ is given by the rotation of
$\pi$ radians around this axis. Of course, the fixed point set of this involution of $H_{n
+ m + 2l - 1}$ consists of $n + m$ arcs. Now, we attach each one of the $m$ extra handles
at two disks centered at the end points of a fixed arc. Finally, we extend the rotation to
the such extra handle as the rotation of $\pi$ radians around its core. Hence, the fixed
arc close up to give a fixed loop.

We remark that $L^{g+1,0}_{g}$ coincides the hyperelliptic involution of $H_g$.

\smallskip

As a consequence of our classification, we see that any orientation preserving involution
of $H_g$ is uniquely determined, up to equivalence, by its restriction to the boundary
$T_g = \Bd H_g$. However, it is worth observing that the restrictions to $T_g$ of
non-equivalent involutions of $H_g$ can be equivalent as involutions of $T_g$, by a PL
homeomorphism of $T_g$ which does not extend to $H_g$. Actually, two involutions of $T_g$
are equivalent if and only if they have the same number of fixed points and they give
raise to quotient surfaces of the same genus $g'$, as it follows from the Hurwitz
classification of branched covering between surfaces (\cite{Hur}, see also \cite{BE}).

\smallskip

From a different point of view, we see that the quotient of $H_g$ under the action of any
orientation preserving involution turns out to be a handlebody $H_{g'}$. Namely, $g' = (g
+ 1)/2 = n + 1$ for $H_{g = 2n + 1} / I_g$ and $g' = (g - n + 1)/2 = m + l$ for $H_{g = n
+ 2m + 2l - 1} / L_g^{n,m}$. Therefore, our result could also be reformulated in terms of
double branched coverings $H_g \to H_{g'}$ between handlebodies.

\section{Preliminaries\label{prelim/sec}}

An {\sl involution} of a PL manifold $X$ is any PL homeomorphism $h: X \to X$ such that
$h \neq \id_X$ and $h^2 = \id_X$. We denote by $\Fix h = \{x \in X \;|\; h(x) = x\}$ the
{\sl fixed point set} of $h$. The involution $h$ is called {\sl free} if $\Fix h =
\emptyset$.

If $h': X' \to X'$ is another involution of the PL manifold $X'$, then we say that $h$ and
$h'$ are {\sl equivalent} if there exists a PL homeomorphism $\eta: X \to X'$ such that
$h' = \eta \circ h \circ \eta^{-1}$.

Here, we focus on orientation preserving involutions the 3-dimensional handlebody $H_g$,
consists of one 0-handle and $g$ orientable 1-handles attached to it, for any $g \geq 0$.
If $h: H_g \to H_g$ is such an orientation preserving involution, then $\Fix h$ is a
(possibly empty) proper PL 1-submanifold of $H_g$. Moreover, the canonical projection
$\pi: H_g \to H_g/h$ turns out to be a double branched covering.

\medskip

In particular, we want to prove the theorem stated in the introduction, providing a
complete classification, up to equivalence, of the orientation preserving involutions of
$H_g$ for any $g \geq 0$.

The proof proceeds by induction on the number $g$ of the 1-handles, starting from the
trivial case of $g = 0$. In this case, we have $H_0 \cong B^3 \subset R^3$, whose only
orientation preserving involution, up to equivalence, is the symmetry $(x,y,z) \mapsto
(x,-y,-z)$ with respect to the $x$-axis (cf. \cite{Liv} and \cite{Wal}), which coincides
with $L_0^{1,0}$.

\medskip

The following lemma concerning involutions of 1-handles, tells us how a given orientation
preserving involution of a disjoint union of orientable handlebodies can be extended to
some extra 1-handles equivariantly attached to it. As an immediate consequence, such
extension is uniquely determined by the equivalence class of the involution induced on the
pairs of attaching disks. This fact will used when performing the inductive step.

\begin{lemma}\label{handle/thm}
The 3-dimensional 1-handle $B^1 \times B^2 \subset R^3$ has only two involutions
preserving the attaching disks $\{-1,1\} \times B^2$, up to equivalence preserving such
disks. Namely, they are the symmetries $(x,y,z) \mapsto (x,-y,-z)$ and $ (x,y,z) \mapsto
(-x,y,-z)$. The first one fixes the core $B^1 \times \{0\}$ of the handle and sends each
attaching disk onto itself, while the second one fixes the diameter $\{0\} \times B^1$ of
the co-core of the handle and swaps the attaching disks.
\end{lemma}

\begin{proof}
Taking into account what we have said about involutions of $B^3$, the lemma can be easily
derived just by considering the possible positions of the arc fixed by the involution with
respect to the attaching disks.
\end{proof}

The other main tools for the inductive step is the next lemma, which allows us to split
any orientation preserving involution of $H_g$ as boundary connected sum of involutions of
simpler handlebodies.

We recall that a properly embedded PL 2-disk $D$ in a bounded 3-manifold $M$ is called
{\sl boundary parallel\/} if there exists a 2-disk $E \subset \Bd M$ such that $\Bd D =
\Bd E$ and $D \cup E$ bounds a 3-cell in $M$. Moreover, if $D'$ is another properly
embedded PL 2-disk in $M$, then $D$ and $D'$ are called {\sl parallel\/} if they are
disjoint and there exists an annulus $A \subset \Bd M$ such that $\Bd A = \Bd D \cup \Bd
D'$ and the 2-sphere $D \cup A \cup D'$ bounds a 3-cell in $M$.

\begin{lemma}\label{split/thm}
Let $h: H_g \to H_g$ be an orientation preserving involution with $g \geq 1$. Then there
exists a properly embedded PL 2-disk $D$ in $H_g$ which is not boundary parallel and such
that either $h(D) \cap D = \emptyset$ or $h(D) = D$ and this meets $\Fix h$ transversally
at one point. In the first case, denoting by $N$ a regular neighborhood of $h(D) \cup D$,
we can assume that $\Cl(H_g - N)$ is PL homeomorphic to $H_{g-2}$ or $H_{g_1} \sqcup
H_{g_2}$ with $g_1 + g_2 = g-1$. In the second case, denoting by $N$ a regular
neighborhood of $h(D) = D$, we have that $\Cl(H_g - N)$ is PL homeomorphic to $H_{g-1}$ or
$H_{g_1} \sqcup H_{g_2}$ with $g_1 + g_2 = g$.
\end{lemma}

\begin{proof}
The first part of the statement follows from Theorem 3 of \cite{GL}. Concerning the second
part, we first observe that $\Cl(H_g - N)$ is a disjoint union of handlebodies (cf.
\cite{Gro}) and $H_g$ can be thought as $\Cl(H_g - N)$ with one (when $h(D) = D$) or two
(when $h(D) \cap D = \emptyset$) 1-handles attached to it. Hence, the only non-trivial
fact to be proved is that $\Cl(H_g - N)$ can be assumed to have at most two components. In
fact, if $h(D) \cap D = \emptyset$ then $\Cl(H_g - N)$ could also have three components,
say $C_1$, $C_2$ and $C_3$. It is not difficult to see that in this case $h$ swaps two of
them, say $C_1$ and $C_2$, and sends the other one onto itself. Since $D$ is not boundary
parallel, $C_1 \cong C_2 \cong H_{g'}$ with $g' \geq 1$. Hence, we can replace the disk
$D$ by a non-separating disk in $C_1$. After that, we have $h(D) \cap D = \emptyset$ and
$\Cl(H_g - N)$ turns out to be connected.
\end{proof}

By previous lemmas, one can easily determine the orientation preserving involutions of
$H_1 \cong S^1 \times B^2 \subset \Bbb C^2$. Since these are known (cf. \cite{Prz} or
\cite{Bar}), we limit ourselves to list them without proof. Up to equivalence, they are
$I_1: (x,y) \mapsto (-x,y)$, $L_1^{0,1}: (x,y) \mapsto (-x,y)$ and $L_1^{2,0}: (x,y)
\mapsto (\widebar x, \widebar y)$, where the bar denotes the complex conjugation, for any
$(x,y) \in S^1 \times B^2$. The first involution is free, while the fixed point sets of
the last two are respectively $S^1 \times \{0\}$ and $\{-1,1\} \times [-1,1]$.

\medskip

We conclude this section by a characterization of the hyperelliptic involutions of $H_g$
for $g \geq 2$. This will be useful in order to simplify the induction argument for
non-free case in the next section.

\begin{lemma}\label{hyperelliptic/thm}
Let $h$ be a non-free orientation preserving involution of $H_g$ with $g \geq 1$. If for
any 2-disk $D$ in $H_g$ given by Lemma \ref{split/thm} the union $h(D) \cup D$ (possibly
coinciding with $D$ itself) disconnect $H_g$, then $h$ is equivalent to $L_g^{g+1,0}$.
\end{lemma}

\begin{proof}
We proceed by induction on $g$. For $g = 0,1$ the statement follows from the above
classification of the orientation preserving involutions of $H_0$ and $H_1$.

Now, assume $g > 1$. Given a disk $D \subset H_g$ as in Lemma \ref{split/thm}, we denote
by $N$ a regular neighborhood of $D \cup h(D)$. Then, $\Cl(H_g - N)$ is disconnected by
hypothesis, and the second part of that lemma implies that $\Cl(H_g - N) = C_1 \sqcup
C_2$, where $C_i \cong H_{g_i}$ for $i=1,2$, with $g_1 + g_2 = g - 1$ if $h(D) = D$ and
$g_1 + g_2 = g - 2$ if $h(D) \cap D = \emptyset$.

Since $h$ is non-free, we have that each of $C_1$ and $C_2$ is sent onto itself by $h$.
Actually, $h$ could in principle swap $C_1$ and $C_2$ (with $g_1 = g_2$) when $h(D) \cap D
= \emptyset$, but in this case it would be free. Moreover, both the restrictions $h_i =
h_{|C_i}: C_i \to C_i$ obviously satisfy the condition of the lemma. Therefore, by the
inductive hypothesis we have $C_i \cong L_{g_i}^{g_i + 1,0}$ for $i = 1,2$.

At this point, can easily conclude that $h \cong L_g^{g+1,0}$ by Lemma \ref{handle/thm},
after observing that $N$ consists of one (if $h(D) = D$) or two (if $h(D) \cap D =
\emptyset$) 1-handles attached to $C_1 \sqcup C_2$ to give $H_g$.
\end{proof}

\section{Proof of the theorem\label{proof/sec}}

Assume first that $h$ is free. Since the Euler characteristic $\chi(H_g) = 1 - g$ is even,
$g = 2n + 1$ for some $n \geq 0$. We will prove that $h \cong I_g$ by induction on $n$,
based on the case $n = 0$, which follows from the above classification of the involutions
of $H_1$.

So, suppose $n > 0$. Let $D \subset H_g$ be a disk as in Lemma \ref{split/thm}. Then $h(D)
\cap D = \emptyset$, being $h$ free. Now, denoting by $N$ a regular neighborhood of $h(D)
\cup D$ and putting $H' = \Cl(H_g - N)$, we have three cases.

\smallskip
{\sl Case 1. $H' \cong H_{g-2}$}. By the inductive hypothesis, $h' = h_{|H'} \cong
I_{g-2}$. Moreover, $N$ consists of a pair of 1-handles equivariantly attached to $H'$,
which are swapped by $h$. Then, up to equivalence, $h$ is the unique possible extension of
$h'$ to $H_g$. Since, up to equivariant PL homeomorphisms, $I_g$ can be obtained in the
same way from $I_{g-2}$, for example by considering as $D$ the leftmost meridian disk in
Figure \ref{involution1/fig}, we have $h \cong I_g$.

\smallskip
{\sl Case 2. $H' = C_1 \sqcup C_2$, with $C_i \cong H_{g_i}$ and $h(C_i) = C_i$ for $i =
1,2$}. Since $g_1 < g$, by the inductive hypothesis $h_{|C_1} \cong I_{g_1}$. Now, if $g_1
> 1$ we know that there exists a disk $D' \subset C_1 \cong H_{g_1}$, such that $C_1 -
(h(D') \cup D')$ is connected. Then, by replacing $D$ with $D'$ thought as a disk in
$H_g$, we are reduced to case 1. On the other hand, if $g_1 = 1$, for any disk $D'$ in
$C_1$, we have that $C_1 - (h(D') \cup D')$ has two components and these are swapped by
$h_{|C_1}$. Then, since also the two attaching disks of the 1-handles given by $N$ are
swapped by $h_{|C_1}$, we can easily conclude that $H_g - (h(D') \cup D')$ is connected.
So, we can once again reduce ourselves to case 1.

\smallskip
{\sl Case 3. $H' = C_1 \sqcup C_2$, with $C_i \cong H_{g_i}$ and $h(C_i) = C_{3-i}$ for $i
= 1,2$}. In this case we have $1 \leq g_1 = g_2 < g$. Then, there exists a disk $D'
\subset C_1$ such that $C_1 - D'$ is connected. Since $h(D') \subset C_2$ and also $C_2 -
h(D')$ is connected (being PL homeomorphic to $C_1 - D'$), we have that $H_g - (h(D') \cup
D')$ is connected too. This allows the reduction to case 1 as above.

\medskip

Now, we assume that $h$ is non-free. We will prove that $h \cong L_g^{n,m}$ by induction
on $g$, based on the cases $g = 0,1$, which follow from the above classification of the
involutions of $H_0$ and $H_1$, and on the cases considered in Lemma
\ref{hyperelliptic/thm}.

So, suppose $g > 1$. Let $D \subset H_g$ be a disk as in Lemma \ref{split/thm}. If for any
such a disk $D$ the union $h(D) \cup D$ disconnect $H_g$, we are done by Lemma
\ref{hyperelliptic/thm}. Hence, we can assume that $H_g - (h(D) \cup D)$ is connected.
Then, denoting by $N$ a regular neighborhood of $h(D) \cup D$ and putting $H' = \Cl(H_g -
N)$, we have $H' \cong H_{g-1}$ if $h(D) = D$ and $H' \cong H_{g-2}$ if $h(D) \cap D =
\emptyset$. We consider this two cases separately.

\smallskip
{\sl Case 1. $h(D) = D$.} By the inductive hypothesis, $h' = h_{|H'} \cong L_{g-1}^{n,m}$
for some $n$ and $m$ such that $1 \leq n + 2m \leq g$. Moreover, $N$ consists of one
1-handle attached to $H'$, whose attaching disks $D_1,D_2 \subset \Bd H'$ are such that
$h'(D_i) = D_i$ and $D_i \cap \Fix h' = \{p_i\} \subset \Int D_i$, for $i = 1,2$. We have
the following two subcases.

\smallskip
{\sl Case 1.1. $p_1$ and $p_2$ are end points of the same arc $A \subset \Fix h'$.} In
this case, when attaching $N$ to $H'$, the arc $A$ closes up to give a fixed loop for $h$.
Now, if $A$ is the rightmost fixed arc in Figure \ref{involution2/fig}, then clearly $h
\cong L_g^{n-1,m+1}$. On the other hand, the half-twists on the disks $E$ and $E' = h'(E)$
in right side of Figure \ref{twists/fig} allows us to equivariantly exchange to
consecutive arcs in $\Fix h'$, hence all the arcs in $\Fix h'$ are equivalent by an
equivariant PL homeomorphisms. Therefore, the final result is the same for any fixed arc
$A \subset \Fix h'$.

\smallskip
{\sl Case 1.2. $p_1$ and $p_2$ are end points of different arcs $A_1, A_2 \subset \Fix
h'$.} In this case, when attaching $N$ to $H'$, the arcs $A_1$ and $A_2$ are joined to
give one fixed arc in $\Fix h$. Now, if $A_1$ and $A_2$ are the rightmost fixed arcs in
Figure \ref{involution2/fig} and the points $p_1$ and $p_2$ are closest the endpoints of
them, then it is not difficult to see that $h \cong L_g^{n-1,m}$. Now, the half-twists on
the disks $E$ and $E' = h'(E)$ in left side of Figure \ref{twists/fig} allows us to
equivariantly exchange the two end points of the same arc in $\Fix h'$. Then, using this
PL homeomorphism, together with that used in the previous case to exchange two consecutive
arcs in $\Fix h'$, we can always equivariantly move the points $p_1$ and $p_2$ in the
preferred position described above. Hence, $h \cong L_g^{n-1,m}$ whatever $p_1$ and $p_2$
are.

\begin{Figure}[htb]{twists/fig}\vskip3pt
\fig{}{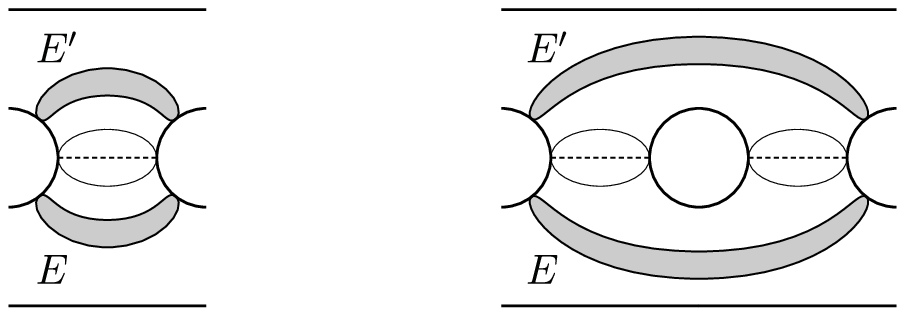}{}{Equivariantly inverting a fixed arc 
       and exchanging two fixed arcs.}\vskip-3pt
\end{Figure}

{\sl Case 2. $h(D) \cap D = \emptyset$.} By the inductive hypothesis, $h' = h_{|H'} \cong
L_{g-2}^{n,m}$ for some $n$ and $m$ such that $1 \leq n + 2m \leq g-1$. Moreover, $N$
consists of a pair of 1-handles equivariantly attached to $H'$, which are swapped by $h$.
Then, up to equivalence, $h$ is the unique possible extension of $h'$ to $H_g$. Since, up
to equivariant PL homeomorphisms, $L_g^{n,m}$ can be obtained in the same way from
$L_{g-2}^{n,m}$, for example by considering as $D$ and $h(D)$ the rightmost meridian disks
in Figure \ref{involution2/fig}, we have $h \cong L_g^{n,m}$.

\end{document}